\documentclass[12pt]{article}
\usepackage{graphicx}
\usepackage{epsfig}
\usepackage{fancyhdr}

\usepackage{amsmath,amsfonts,amssymb,amsthm}
\theoremstyle{plain}
\newtheorem{thm}{Theorem}
\newtheorem*{thm*}{Theorem}

\newtheorem*{lem*}{Lemma}

\theoremstyle{definition}

\newcommand{\reftit}{\textit}    
\newcommand{\refis}{\textbf}     

\begin{document}

\title{Orthogonality and probability: beyond nearest neighbor transitions}
\author{Yevgeniy Kovchegov\footnote{
 Department of Mathematics,  Oregon State University, Corvallis, OR  97331-4605, USA
 \texttt{kovchegy@math.oregonstate.edu}}}
\date{ }
\maketitle

\abstract{In this article, we will explore why Karlin-McGregor method of using orthogonal polynomials in the study of Markov
 processes was so successful for one dimensional nearest neighbor processes, but failed beyond nearest neighbor transitions.
 We will proceed by suggesting and testing  possible fixtures.}

\section{Introduction}

 This paper was influenced by the approaches described in Deift \cite{deift} and questions considered in Gr\"unbaum \cite{g}.
 
 The Karlin-McGreogor diagonalization can be used to answer recurrence/transience questions, as well as those of
  probability harmonic functions, occupation times and hitting times, and a large number of other quantities obtained by solving
 various recurrence relations, in the study of Markov chains,   see \cite{km1}, \cite{km2}, \cite{km2a},  \cite{km3}, \cite{karlin}, \cite{szego}, \cite{schoutens}, \cite{dksc}, \cite{kmn}.  However with some   exceptions (see \cite{km4}) those were nearest neighbor Markov chains on half-line. Gr\"unbaum \cite{g} mentions two main drawbacks to the method as (a) ``typically one cannot get either the polynomials or the measure explicitly", and (b) ``the method is restricted
   to `nearest neighbour' transition probability chains that give rise to tridiagonal matrices and thus to orthogonal polynomials".
   In this paper we attempt to give possible answers to the second question of Gr\"unbaum \cite{g} for general reversible Markov
   chains. In addition, we will consider possible applications of the newer methods in orthogonal polynomials such as using 
   Riemann-Hilbert approach, see \cite{deift}, \cite{deift1} and \cite{ku}, and their probabilistic interpretations.

 In Section 2, we will give an overview of the Karlin-McGregor method from a naive college linear algebra perspective. In 2.3, we will
 give a Markov chain interpretation to the result of Fokas , Its and Kitaev, connecting orthogonal polynomials and Riemann-Hilbert problems. Section 3 deals with one dimensional random walks with jumps of size $\leq m$, the $2m+1$ diagonal operators.
 There we consider dioganalizing with orthogonal functions. In 3.2, as an example we consider a pentadiagonal operator and use Plemelj 
 formula,  and a two sided interval to obtain the respective diagonalization. In Section 4, we use the constructive approach of 
 Deift \cite{deift} to produce the Karlin-McGregor diagonalization for all irreducible reversible Markov chains. After that, we revisit the example from Section 3.

\section{Eigenvectors of probability operators} 

Suppose $P$ is a tridiagonal operator of a one-dimensional Markov chain on $\{0,1,\dots\}$ with forward probabilities $p_k$
 and backward probabilities $q_k$. Suppose $\lambda$ is an eigenvalue of $P$ and 
${\bf q}^T(\lambda)=\left(\begin{array}{c}Q_0 \\ Q_1 \\ Q_2 \\ \vdots \end{array}\right)$ is the corresponding right eigenvector such that $Q_0=1$.
So $\lambda {\bf q}^T=P{\bf q}^T$ generates the recurrence relation for $Q_j$. Then each $Q_j(\lambda)$ is a polynomial of $j$-th degree.
The Karlin-McGregor method derives the existance of a probability distribution $\psi$ such that polynomials $Q_j(\lambda)$ are
orthogonal with respect to $\psi$. In other words, if $\pi$ is stationary with $\pi_0=1$ and $<\cdot,\cdot>_{\psi}$ is the inner product
 in $L^2(d\psi)$, then
$$<Q_i,Q_j>_{\psi}={\delta_{i,j} \over \pi_j}$$
Thus $\{\sqrt{\pi_j} Q_j(\lambda)\}_{j=0,1,\dots}$ are orthonormal polynomials, where $\pi_j={p_0 \dots p_{j-1} \over q_1 \dots q_j}$. 
Also observe from the recurrence relation that the leading coefficient of $Q_j$ is ${1 \over p_0\dots p_{j-1}}$.
\vskip 0.2 in
\noindent
Now, $\lambda^t {\bf q}^T=P^t{\bf q}^T$ implies $\lambda^t q_i =(P^t {\bf q}^T)_i$ for each $i$, and
$$<\lambda^t Q_i, Q_j>_{\psi}=<(P^t {\bf q}^T)_i,Q_j>_{\psi}={p_t(i,j) \over \pi_j}$$
Therefore
$$p_t(i,j)=\pi_j <\lambda^t Q_i, Q_j>_{\psi}$$

 Since the spectrum of $P$ lies entirely inside $[-1,1]$ interval, then so is the support of $\psi$. Hence, for $|z|>1$,
 the generating function 
 $$G_{i,j}(z)=\sum_{t=0}^{+\infty} z^{-t} p_t(i,j)=-z \pi_j <{Q_i \over \lambda-z}, Q_j>_{\psi}
 =-z \pi_j  \int {Q_i(\lambda) Q_j(\lambda) \over \lambda-z} d\psi(\lambda)$$

\subsection{Converting to a Jacobi operator}
  Let $b_k=\sqrt{\pi_k \over \pi_{k+1}}p_k$, then $b_k=\sqrt{\pi_{k+1} \over \pi_k}q_{k+1}$ due to reversibility condition. Thus the recurrence relation for ${\bf q}$,
  $$\lambda \sqrt{\pi_k} Q_k=q_k \sqrt{\pi_k} Q_{k-1}+(1-q_k-p_k) \sqrt{\pi_k} Q_k+p_k  \sqrt{\pi_k} Q_{k+1}~,$$
 can be rewritten as
  $$\lambda \sqrt{\pi_k} Q_k=b_{k-1} \sqrt{\pi_k} Q_{k-1}+a_k \sqrt{\pi_k} Q_k+b_k  \sqrt{\pi_k} Q_{k+1},$$
 where $a_k= 1-q_k-p_k$.
  Therefore ${\bf \widetilde{q}}=(\sqrt{\pi_0} Q_0, \sqrt{\pi_1} Q_1, \dots )$ solves $\widetilde{P} {\bf \widetilde{q}} = \lambda {\bf \widetilde{q}}$, where
  $$\widetilde{P}=\left(\begin{array}{cccc}a_0 & b_0 & 0 & \dots \\b_0 & a_1 & b_1 & \ddots \\0 & b_1 & a_2 & \ddots \\\vdots & \ddots & \ddots & \ddots\end{array}\right)$$
  is a Jacoby (symmetric triangular with $b_k>0$) operator. Observe that $\widetilde{P}$ is self-adjoint.

 The above approach extends to all reversible Markov chains. Thus every reversible Markov operator is equivalent to
 a self-adjoint operator, and therefore has an all real spectrum.

\subsection{Karlin-McGregor: a simple picture}
 It is a basic fact from linear algebra that if $\lambda_1, \dots, \lambda_n$ are distinct real eigenvalues of
 an $n \times n$ matrix $A$, and if $u_1,\dots,u_n$ and $v_1,\dots, v_n$ are  the corresponding left and right eigenvectors.
 Then $A$ diagonalizes  as follows
 $$A^t=\sum_{j} {\lambda^t v^T_j u_j \over u_j v^T_j}=\int_{\sigma(A)} \lambda^t v^T(\lambda) u(\lambda) d \psi(\lambda)~,$$
 where $u(\lambda_j)=u_j$,  $v(\lambda_j)=v_j$, spectrum $\sigma(A)=\{\lambda_1, \dots, \lambda_n \}$, and
 $$\psi(\lambda)=\sum_{j} { 1 \over u(\lambda) v^T(\lambda)} \delta_{\lambda_j} (\lambda)={n \over u(\lambda) v^T(\lambda)} U_{\sigma(A)}(\lambda)$$
Here $U_{\sigma(A)}(\lambda)$ is the uniform distribution over the spectrum $\sigma(A)$.

 It is {\bf important} to observe that the above integral representation is only possible if $u(\lambda)$ and $v(\lambda)$ are well
 defined - each eigenvalue has multiplicity one, i.e. all distinct real eigenvalues. As we will see later, this will become crucial 
 for Karlin-McGregor diagonalization of reversible Markov chains. The operator for a reversible Markov chain is bounded and 
 is equivalent to a  self-adjoint operator,  and as such has a real bounded spectrum. However the eigenvalue multiplicity will determine whether the operator's diagonalization can be expressed in a form of a spectral integral.

 Since the spectrums $\sigma(P)=\sigma(P^*)$, we will extend the above diagonalization identity to the operator $P$ in the separable Hilbert space $l^2(\mathbb{R})$.  First, observe that ${\bf u}(\lambda)=(\pi_0 Q_0, \pi_1 Q_1, \dots)$ satisfies
 $${\bf u}P =\lambda P$$
 due to reversibility. Hence, extending from a finite case to an infinite dimensional space $l^2(\mathbb{R})$,  obtain
 $$P^t=\int \lambda^t {\bf q}^T(\lambda) {\bf u}(\lambda) d \psi(\lambda)
 =\int \lambda^t \left(\begin{array}{ccc}\pi_0 Q_0 Q_0 & \pi_1 Q_0 Q_1 & \cdots \\\pi_0 Q_1 Q_0 & \pi_1 Q_1 Q_1 & \cdots \\\vdots & \vdots & \ddots\end{array}\right) d \psi(\lambda)~,$$  
 where
 $$\psi(\lambda)= \lim_{n \rightarrow +\infty} \psi_n(\lambda)$$
 The above is the weak limit of
 $$\psi_n(\lambda)={n \over {\bf u}(\lambda) {\bf q}^T(\lambda)} U_{\sigma(A_n)}(\lambda)~,$$
 where $A_n$ is the restriction of $P$ to the first $n$ coordinates, $<e_0,\dots,e_{n-1}>$
 $$A_n= \left(\begin{array}{ccccc}1-p_0 & p_0 & 0 & \cdots & 0 \\q_1 & 1-q_1-p_1 & p_1 & \ddots & \vdots \\0 & q_2 & 1-q_2 & \ddots & 0 \\\vdots & \ddots & \ddots & \ddots & p_{n-2} \\0 & \cdots & 0 & q_{n-1} & 1-q_{n-1}-p_{n-1}\end{array}\right)$$
  Observe that if $Q_n(\lambda)=0$ then $(Q_0(\lambda),\dots,Q_{n-1}(\lambda))^T$ is the corresponding right eigenvector of $A_n$. 
 Thus the spectrum of $\sigma(A_n)$ is the roots of $$Q_n(\lambda)=0$$ 
 So 
 $$\psi_n(\lambda)={n \over {\bf u}(\lambda) {\bf q}^T(\lambda)} U_{Q_n=0}(\lambda)
 ={n \over \sum_{k=0}^{n-1} \pi_k Q_k^2(\lambda)} U_{Q_n=0}(\lambda)~.$$ 
 The orthogonality follows if we plug in $t=0$. Since $\pi_0 Q_0 Q_0 =1$, $\psi$ should integrate to one.

 \vskip 0.3 in
 \noindent
 {\bf Example.} {\it Simple random walk and Chebyshev polynomials.} The Chebyshev polynomials of the first kind are the ones characterizing a one dimensional simple random walk on half line, i.e. the ones with generator 
$$P_{ch}=\left(\begin{array}{ccccc}0 & 1 & 0 & 0 & \cdots \\{1 \over 2} & 0 & {1 \over 2} & 0 & \cdots \\0 & {1 \over 2} & 0 & {1 \over 2} & \ddots \\0 & 0 & {1 \over 2} & 0 & \ddots \\\vdots & \vdots & \ddots & \ddots & \ddots\end{array}\right)$$
So, $T_0(\lambda)=1$, $T_1(\lambda)=\lambda$ and $T_{k+1}(\lambda)=2 \lambda T_k(\lambda)-T_{k-1}(\lambda)$ for $k=2,3, \dots$.
 The Chebyshev polynomials satisfy the following trigonometric identity:
 $$T_k(\lambda)=\cos(k \cos^{-1}(\lambda))$$  
 Now, $$\psi_n(\lambda)={n \over \sum_{k=0}^{n-1} \pi_k T_k^2(\lambda)} U_{\{\cos(n \cos^{-1}(\lambda))=0\}}(\lambda)~,$$
 where $\pi(0)=1$ and $\pi(1)=\pi(2)=\dots=2$.
 Here $$U_{\{\cos(n \cos^{-1}(\lambda))=0\}}(\lambda)=U_{\{\cos^{-1}(\lambda))={\pi \over 2n}+{\pi k \over n}, ~k=0,1,\dots,n-1\}}(\lambda)$$
 Thus if $X_n \sim U_{\{\cos(n \cos^{-1}(\lambda))=0\}}$, then $Y_n=\cos^{-1}(X_n) \sim U_{\{{\pi \over 2n}+{\pi k \over n}, ~k=0,1,\dots,n-1\}}$ 
 and $Y_n$ converges weakly to $Y \sim U_{[0,\pi]}$. Hence $X_n$ converges weakly to
 $$X=\cos(Y) \sim {1 \over \pi \sqrt{1 -\lambda^2}} \chi_{[-1,1]}(\lambda) d\lambda~,$$
 i.e. $$U_{\{\cos(n \cos^{-1}(\lambda))=0\}}(\lambda) \rightarrow {1 \over \pi \sqrt{1 -\lambda^2}} \chi_{[-1,1]}(\lambda) d\lambda$$
 Also observe that if $x=\cos(\lambda)$, then
 $$ \sum_{k=0}^{n-1} \pi_k T_k^2(\lambda)=-1+2\sum_{k=0}^{n-1} \cos^2(kx)=n-{1 \over 2}+{\sin((2n-1)x) \over 2\sin(x)} $$
 Thus
 $$d\psi_n(\lambda) \rightarrow d \psi(\lambda) ={1 \over \pi \sqrt{1 -\lambda^2}} \chi_{[-1,1]}(\lambda) d\lambda$$

\subsection{Riemann-Hilbert problem and a generating function of $p_t(i,j)$}
Let us write $\sqrt{\pi_j} Q_j(\lambda)=k_j P_j(\lambda)$, where $k_j={1 \over \sqrt{p_0 \dots p_{j-1}} \sqrt{q_1 \dots q_j}}$
 is the leading coefficient of $\sqrt{\pi_j} Q_j(\lambda)$, and $P_j(\lambda)$ is therefore a {\it monic} polynomial.
 
 In preparation for the next step, let $w(\lambda)$ be the probability density function associated with the spectral measure $\psi$: 
 $d \psi(\lambda)=w(\lambda) d \lambda$ on the compact support, $supp(\psi) \subset [-1,1]=\Sigma$. Also let
 $$C(f)(z)={1 \over 2\pi i}\int_{\Sigma} {f(\lambda) \over \lambda -z} d \psi(\lambda)$$ 
 denote the Cauchy transform w.r.t. measure $\psi$.
 
  First let us quote the following theorem.
  \begin{thm*} {\bf[Fokas, Its and Kitaev, 1990]} Let $$v(x)=\left(\begin{array}{cc}1 & w(x) \\0 & 1\end{array}\right)$$
  be the jump matrix. Then, for any $n \in \{0,1,2,\dots \}$,
  $$m^{(n)}(z)=\left(\begin{array}{cc}P_n(z) & C(P_n w)(z) \\-2\pi i k^2_{n-1}P_{n-1}(z)  & -2\pi i k^2_{n-1}C(P_{n-1} w)(z) \end{array}\right),~\text{ for all } z \in  \mathbb{C} \setminus \Sigma,$$
 is the unique solution to the Riemann-Hilbert problem with the above jump matrix $v(x)$ and $\Sigma$ that satisfies
  the following  condition
  \begin{equation} \label{RHcondition}
  m^{(n)}(z) \left(\begin{array}{cc}z^{-n} & 0 \\0 & z^n\end{array}\right) \rightarrow I~\text{ as }~z \rightarrow \infty~.
  \end{equation}
  \end{thm*}

 The Riemann-Hilbert problem, for an oriented smooth curve $\Sigma$, is the problem of finding $m(z)$, analytic in 
 $\mathbb{C} \setminus \Sigma$ such that
 $$m_+(z)=m_-(z)v(z),~~~\text{ for all }z \in \Sigma,$$
 where $m_+$ and $m_-$ denote respectively the limit from the left and the limit from the right, for the function $m$, as we
 approach a point on $\Sigma$. 
 
 Suppose we are given the weight function $w(\lambda)$ for the Karlin-McGregor orthogonal polynomials ${\bf q}$.
 If $m^{(n)}(z)$ is the solution of the Riemann-Hilbert problem as in the above theorem, then for $|z|>1$,
 $$m^{(n)}(z)=\left(\begin{array}{cc}{1 \over k_n\sqrt{\pi_n}} Q_n(z) & -{1 \over 2\pi i k_n \sqrt{\pi_n} z^{n+1}}G_{0,n} \\-2\pi i{ k_{n-1} \over \sqrt{\pi_{n-1}}}Q_{n-1}(z)  & { k_{n-1} \over \sqrt{\pi_{n-1}}z^n} G_{0,n-1} (z) \end{array}\right)$$
 $$=\left(\begin{array}{cc}q_1\dots q_nQ_n(z) & -{q_1\dots q_n \over 2\pi i z^{n+1}}G_{0,n} \\{-2\pi i \over p_0\dots p_{n-2}}Q_{n-1}(z)  & { 1 \over p_0\dots p_{n-2} z^n} G_{0,n-1} (z) \end{array}\right)$$

\section{Beyond nearest neighbor transitions}

 Observe that the Chebyshev polynomials were used to diagonalize a simple one dimensional random walk reflecting
 at the origin. Let us consider a random walk where jumps of sizes one and two are equiprobable
 $$P=\left(\begin{array}{ccccccc}0 & {1 \over 2} & {1 \over 2} & 0 & 0 & 0 & \dots \\{1 \over 4} & {1 \over 4} & {1 \over 4} & {1 \over 4} & 0 & 0 & \dots \\{1 \over 4} & {1 \over 4} & 0 & {1 \over 4} & {1 \over 4} & 0 & \dots \\0 & {1 \over 4} & {1 \over 4} & 0 & {1 \over 4} & {1 \over 4} & \ddots \\0 & 0 & {1 \over 4} & {1 \over 4} & 0 & {1 \over 4} & \ddots \\0 & 0 & 0 & {1 \over 4} & {1 \over 4} & 0 & \ddots \\\cdots & \cdots & \cdots & \ddots & \ddots & \ddots & \ddots\end{array}\right)$$
  The above random walk with the reflector at the origin is reversible with $\pi(0)=1$ and $\pi(1)=\pi(2)=\dots=2$.
  The Karlin-McGregor representation with orthogonal polynomials will not automatically extend to this case. However this does not rule out obtaining a Karlin-McGregor diagonalization with orthogonal functions.

 In the case of the above pentadiagonal Chebyshev operator, some eigenvalues will be of geometric multiplicity two as
 $$P=P^2_{ch}+{1 \over 2}P_{ch}-{1 \over 2}I~,$$
 where $P_{ch}$ is the original tridiagonal Chebyshev operator.

\subsection{$2m+1$ diagonal operators}  

 Consider a $2m+1$ diagonal reversible probability operator $P$. Suppose it is
 Karlin-McGregor diagonalizable. Then for a given $\lambda \in \sigma(P)$,  let 
 ${\bf q}^T(\lambda)=\left(\begin{array}{c}Q_0 \\ Q_1 \\ Q_2 \\ \vdots \end{array}\right)$ once again denote the corresponding right eigenvector such that $Q_0=1$. Since the operator is more than tridiagonal, we encounter the problem of finding the next $m-1$ functions,
 $Q_1(\lambda)=\mu_1(\lambda)$, $Q_2(\lambda)=\mu_2(\lambda), \dots, Q_{m-1}(\lambda)=\mu_{m-1}(\lambda)$. 
 
 Observe that  ${\bf q}={\bf q_0}+{\bf q_1} \mu_1+\dots+{\bf q_{m-1}} \mu_{m-1}$, where each
 ${\bf q}_j^T(\lambda)=\left(\begin{array}{c}Q_{0,j} \\ Q_{1,j} \\ Q_{2,j} \\ \vdots \end{array}\right)$
 solves $P{\bf q}_j^T=\lambda {\bf q}_j^T$ recurrence relation with the initial conditions 
 $$Q_{0,j}(\lambda)=0,~~ \dots,~~ Q_{j-1,j}(\lambda)=0,~~ Q_{j,j}(\lambda)=1,~~ Q_{j+1,j}(\lambda)=0,~~ \dots,~~ Q_{m-1,j}(\lambda)=0$$
 In other words,
 ${\bf q}^T(\lambda)={\bf Q}(\lambda) \mu^T~,$
  where 
  ${\bf Q}(\lambda)=\left[\begin{array}{cccc} &  &  &  \\| & | &  & | \\{\bf q}_0^T & {\bf q}_1^T & \cdots & {\bf q}_{m-1}^T \\| & | &  & | \\ &  &  & \end{array}\right]$
 and $\mu^T=\left(\begin{array}{c}1 \\\mu_1(\lambda) \\\vdots \\\mu_{m-1}(\lambda)\end{array}\right)$ is such that
 ${\bf q}(\lambda) \in l^2(\mathbb{R})$ for each $\lambda \in \sigma(P)$.
 
  Let again $A_n$ denote the restriction of $P$ to the first $n$ coordinates, $<e_0,\dots,e_{n-1}>$
  Observe that if $Q_n(\lambda)=\dots=Q_{n+m-1}(\lambda)=0$ then $(Q_0(\lambda),\dots,Q_{n-1}(\lambda))^T$ is the corresponding right eigenvector of $A_n$. 
 Thus the spectrum of $\sigma(A_n)$ consists of the roots of 
 $$\det\left(\begin{array}{cccc}Q_{n,0}(\lambda) & Q_{n,1}(\lambda) &  & Q_{n,m-1}(\lambda) \\Q_{n+1,0}(\lambda) & Q_{n+1,1}(\lambda) &  & Q_{n+1,m-1}(\lambda) \\\vdots & \vdots & \cdots & \vdots \\Q_{n+m-1,0}(\lambda) & Q_{n+m-1,1}(\lambda) &  & Q_{n+m-1,m-1}(\lambda)\end{array}\right)=0$$

\subsection{Chebyshev operators}

 Let us now return to the example generalizing the simple random walk reflecting at the origin.  There one step and two step 
 jumps were equally likely. The characteristic equation $z^4+z^3-4\lambda z^3+z^2+z=0$ for the recurrence relation
 $$c_{n+2}+c_{n+1}-4\lambda c_n +c_{n-1}+c_{n-2}=0$$
 can be easily solved by observing that if $z$ is a solution then so are $\bar{z}$ and ${1 \over z}$. 
 The solution in radicals is expressed as
 $~z_{1,2}=r_1 \pm i \sqrt{1-r_1^2}~$ and $~z_{3,4}=d_2 \pm i \sqrt{1-r_2^2}~$,
 where
 $r_1={-1 + \sqrt{9+16\lambda} \over 4}$ and 
$ r_2={-1 - \sqrt{9+16\lambda} \over 4}$.

Observe that $r_1$ and $r_2$ are the two roots of $s(x)=\lambda$, where $s(x)=x^2+{1 \over 2}x-{1 \over 2}$ is the
polynomial for which $$P=s(P_{ch})$$  
In general, the following is true for all operators $P$ that represent symmetric random walks reflecting at the origin, and that allow jumps of up to $m$ flights: there is a polynomial $s(x)$ such that $P=s(P_{ch})$ and the roots $z_j$ of the characteristic relation in
$\lambda {\bf c} = P {\bf c}$ will lie on a unit circle with their real parts $Re(z_j)$ solving $s(x)=\lambda$. The reason for the latter is
the symmetry of the corresponding characteristic equation of order $2m$, implying ${1 \over z_j}=\bar{z_j}$, and therefore the characteristic equation for $\lambda {\bf c} = P {\bf c}$ can be rewritten as
$$s\left( {1 \over 2}\left[z+{1 \over z} \right]\right)=\lambda~,$$  
where ${1 \over 2}\left[z+{1 \over z} \right]$ is the Zhukovskiy function.

 In our case, $s(x)=\left(x+{1 \over 4}\right)^2-{9 \over 16}$, and for $\lambda \in \left(-{9 \over 16},0 \right]$, there will be two candidates 
 for $\mu_1(\lambda)$,
 $$\mu_+(\lambda)=r_1={-1 + \sqrt{9+16\lambda} \over 4}~~~\text{ and }~~\mu_-(\lambda)=r_2={-1 - \sqrt{9+16\lambda} \over 4}$$
 Taking $0 \leq \arg{z}<2\pi$ branch of the logarithm $\log{z}$, and applying Plemelj formula, one would obtain 
 $$\mu_1(z)=-{1 \over 4}+z^{1 \over 2} \exp\left\{{1 \over 2} \int_{-{9 \over 16}}^0 {ds \over s-z}\right\}~,$$
 where $\mu_+(\lambda)=\lim_{z \rightarrow \lambda,~Im(z)>0} \mu_1(z)$ and 
 $\mu_-(\lambda)=\lim_{z \rightarrow \lambda,~Im(z)<0} \mu_1(z)$.
 
Now, as we defined $\mu_1(z)$, we can propose the limits of integration to be a contour in $\mathbb{C}$ 
consisting of  
$\left[-{9 \over 16},0\right)_+=\lim_{\varepsilon \downarrow 0} \left\{z=x+i\varepsilon~:~x \in \left[-{9 \over 16},0\right) \right\}$,
and $\left[-{9 \over 16},0\right)_-=\lim_{\varepsilon \downarrow 0} \left\{z=x-i\varepsilon~:~x \in \left[-{9 \over 16},0\right) \right\}$,
and the $[0,1]$ segment. Then
$$P^t=\int_{\left[-{9 \over 16},0\right)_- \cup \left[-{9 \over 16},0\right)_+ \cup [0.1]} \lambda^t {\bf q}^T(\lambda) {\bf u}(\lambda) d \psi(\lambda),$$
where ${\bf u}(\lambda)$ is defined as before, and
 $$d\psi(\lambda)={1 \over 2\pi \sqrt{\lambda+{9 \over 16}}}\left( {\chi_{[-{9 \over 16},0)_-}(\lambda) \over \sqrt{1-\left(\sqrt{\lambda+{9 \over 16}}+{1 \over 4}\right)^2}} 
 +{\chi_{[-{9 \over 16},0)_+}(\lambda) +\chi_{[0,1]}(\lambda) \over \sqrt{1-\left(\sqrt{\lambda+{9 \over 16}}-{1 \over 4}\right)^2}}\right)d\lambda$$

  Let us summarize this section as follows. If the structure of the spectrum does not allow Karlin-McGregor diagonalization with 
  orthogonal functions over $[-1,1]$, say when there are two values of $\mu^T(\lambda)$ for some $\lambda$, then one may use 
  Plemelj formula to obtain an integral diagonalization of $P$ over the corresponding two sided interval.

\section{Spectral Theorem and why orthogonal polynomials work}

 The constructive proofs in the second chapter of Deift \cite{deift} suggest the reason why Karlin-McGregor theory of diagonalizing 
 with orthogonal  polynomials works for all time reversible Markov chains. Using the same logical steps as in \cite{deift}, we can 
 construct a map ${\cal M}$ which assigns a probability measure $d\psi$ to a reversible transition operator $P$ on a countable state space 
 $\{0,1,2,\dots\}$.
 W.l.o.g. we can assume $P$ is symmetric as one can instead consider
 $$\left(\begin{array}{ccc}\sqrt{\pi_0} & 0 & \cdots \\0 & \sqrt{\pi_1} & \ddots \\\vdots & \ddots & \ddots\end{array}\right)
 P\left(\begin{array}{ccc}{1 \over \sqrt{\pi_0}} & 0 & \cdots \\0 & {1 \over \sqrt{\pi_1}} & \ddots \\\vdots & \ddots & \ddots\end{array}\right)$$
 which is symmetric, and its spectrum coinciding with spectrum $\sigma(P) \subset [-1,1]$.

 Now, for $z \in \mathbb{C} \setminus \mathbb{R}$ let $G(z)=(e_0, (P-zI)^{-1}e_0)$. Then
 $$ImG(z)={1 \over 2i}\left[(e_0, (P-zI)^{-1}e_0)-(e_0, (P-\bar{z}I)^{-1}e_0)\right]=(Im(z))|(P-zI)^{-1}e_0|^2$$
  and therefore $G(z)$ is a {\bf Herglotz function}, i.e. $G(z)$ is an analytic map from $\{Im(z)>0\}$ into $\{Im(z)>0\}$,
  and as all such functions, it can be represented as
  $$G(z)=az+b+\int_{-\infty}^{+\infty} \left({1 \over s-z} -{s \over s^2+1} \right) d\psi(s),~~Im(z)>0$$
  In the above representation $a \geq 0$ and $b$ are real constants and $d\psi$ is a Borel measure such that
  $$\int_{-\infty}^{+\infty} {1 \over s^2+1} d\psi(s) < \infty$$
  Deift \cite{deift} uses $G(z)=(e_0, (P-zI)^{-1}e_0)=-{1 \over z}+O(z^{-2})$ to show $a=0$ in our case, and 
  $$b=\int_{-\infty}^{+\infty} {s \over s^2+1} d\psi(s)$$
  as well as the uniqueness of $d\psi$. Hence
  $$G(z)=\int{d\psi(s) \over s-z},~~Im(z)>0$$

 The point of all these is to construct the spectral map 
 $$\mathcal{M}: \{\text{reversible Markov operators P} \} \rightarrow \{\text{probability measures }\psi\text{ on }[-1,1]\text{ with compact }supp(\psi)\}$$
 The asymptotic evaluation of both sides in
 $$(e_0, (P-zI)^{-1}e_0)=\int{d\psi(s) \over s-z},~~Im(z)>0$$
  implies
  $$(e_0, P^k e_0)=\int s^k d\psi(s)$$
  
 Until now we were reapplying the logical steps in Deift \cite{deift} for the case of reversible Markov chains. 
 However, in the original, the second chapter of Deift \cite{deift} gives a constructive proof of  the following spectral
 theorem, that summarizes as 
 $$\mathcal{U}: \{\text{bounded Jacobi operators on }l^2 \} \rightleftharpoons \{\text{probability measures }\psi\text{ on }\mathbb{R}\text{ with compact }supp(\psi)\},$$
where $\mathcal{U}$ is one-to-one onto.

\begin{thm}\label{spectralJacobi}{\bf [Spectral Theorem]}
For every bounded Jacobi operator $\mathcal{A}$ there exists a unique probability measure $\psi$ with compact support such that
$$G(z)=\left(e_0,  (\mathcal{A}-zI)^{-1}e_0 \right)=\int_{-\infty}^{+\infty} {d \psi(x) \over x-z}$$
The spectral map  $\mathcal{U}:\mathcal{A} \rightarrow d\psi$ is one-to-one onto, and for every $f \in L^2(d\psi)$,
$$(\mathcal{UAU}^{-1} f)(s)=sf(s)$$ in the following sense
$$(e_0, \mathcal{A} f(\mathcal{A})e_0)=\int s f(s) d\psi(s)$$
\end{thm}

 \noindent
 So suppose $P$ is a reversible Markov chain, then 
 $$\mathcal{M}:P \rightarrow d\psi~~~\text{and}~~~\mathcal{U}^{-1}: d\psi \rightarrow P_{\triangle}~,$$
 where $P_{\triangle}$ is a unique Jacobi operator such that
 $$(e_0,P^k e_0)=\int s^k d\psi(s)=(e_0,P_{\triangle}^k e_0)$$
 Now, if $Q_j(\lambda)$ are the orthogonal polynomials w.r.t. $d\psi$ associated with $P_{\triangle}$, 
 then $Q_j(P_{\triangle})e_0=e_j$ and
 $$\delta_{i,j}=(e_i,e_j)=(Q_i(P_{\triangle})e_0,Q_j(P_{\triangle})e_0)=(Q_i(P)e_0,Q_j(P)e_0)$$
 Thus, if $P$ is irreducible, then  $f_j=Q_j(P)e_0$ is an orthonormal basis for Karlin-McGregor diagonalization.
 If we let $F=\left[\begin{array}{ccc}| & | &  \\ f^T_0 & f^T_1 & \cdots \\| & | & \end{array}\right]$, then
 $$P^t=\left(\begin{array}{ccc}  &   &   \\  & (P^t e_i,e_j) &   \\  &   &  \end{array}\right)=F\left(\begin{array}{ccc}  &   &   \\  & \int_{-1}^1 s^t Q_i(s)Q_j(s) d\psi(s) &   \\  &   &  \end{array}\right)F^{T}, $$
 where $F^T=F^{-1}$.
 Also Deift \cite{deift} provides a way for constructing
 $$\mathcal{U}^{-1} \mathcal{M}: P \rightarrow P_{\triangle}$$ 
 Since $P_{\triangle}$ is a Jacobi operator, it can be represented as
 $$P_{\triangle}=\left(\begin{array}{cccc}a_0 & b_0 & 0 & \cdots \\b_0 & a_1 & b_1 & \ddots \\0 & b_1 & a_2 & \ddots \\\vdots & \ddots & \ddots & \ddots\end{array}\right)~~~~b_j>0$$
 Now,
 $$(e_0,Pe_0)=(e_0,P_{\triangle}e_0)=a_0,~~~(e_0,P^2e_0)=(e_0,P_{\triangle}^2e_0)=a_0^2+b_0^2$$
 $$(e_0,P^3e_0)=(e_0,P^3_{\triangle}e_0)=(a_0^2+b_0^2)a_0+(a_0+a_1)b_0^2$$
 $$\text{ and }(e_0,P^4e_0)=(e_0,P^4_{\triangle}e_0)=(a_0^2+b_0^2)^2+(a_0+a_1)^2b_0^2+b_0^2b_1^2$$
 thus providing us with the coefficients of the Jacobi operator, $a_0$, $b_0$, $a_1, \dots$, and therefore the orthogonal
 polynomials $Q_j$.
 
  \vskip 0.3 in
 \noindent
 {\bf Example.} {\it Pentadiagonal Chebyshev operator.}
 For the pentadiagonal $P$ that represents the symmetric random walk with equiprobable jumps of sizes one and two,
 $$(e_0,Pe_0)=0,~~(e_0,P^2e_0)={1 \over 4},~~(e_0,P^3e_0)={3 \over 32},~~(e_0,P^4e_0)={9 \over 64},~~ \dots$$
 Thus $$a_0=0,~~b_0={1 \over 2},~~a_1={3 \over 8},~~b_1={\sqrt{11} \over 8}, \text{ etc. }$$
 So
 $$P_{\triangle}=\left(\begin{array}{cccc} 0 & {1 \over 2} & 0 & \cdots \\{1 \over 2} & {3 \over 8} & {\sqrt{11} \over 8} & \ddots \\0 & {\sqrt{11} \over 8} & \ddots & \ddots \\\vdots & \ddots & \ddots & \ddots\end{array}\right)$$
 and
 $$Q_0(\lambda)=1,~~~Q_1(\lambda)=2 \lambda,~~~
 Q_2(\lambda)={32 \over \sqrt{11}} \lambda^2-{6 \over \sqrt{11}}\lambda-{4 \over \sqrt{11}}, ~~\dots$$
 Then applying classical Fourier analysis, one would obtain
 $$\left(e_0,  (P-zI)^{-1}e_0 \right)={1 \over 2\pi} \int_0^{2\pi} {d\theta \over {1 \over 2}[\cos(\theta)+\cos(2\theta)]-z}
 =\int_{-{9 \over 16}}^1 {d\psi(s) \over s-z}~,$$
 where
 $$d\psi(s)={1 \over 2\pi \sqrt{s+{9 \over 16}}}\left( {\chi_{[-{9 \over 16},1]}(s) \over \sqrt{1-\left(\sqrt{s+{9 \over 16}}-{1 \over 4}\right)^2}} 
 +{\chi_{[-{9 \over 16},0)}(s) \over \sqrt{1-\left(\sqrt{s+{9 \over 16}}+{1 \over 4}\right)^2}}\right)ds$$
  To obtain the above expression for $d\psi$ we used the fact that $\left(e_0,  (P-zI)^{-1}e_0 \right)$ would be the same if there were no reflector at zero.

 \subsection{Applications of Karlin-McGregor diagonalization}
 
 Let us list some of the possible applications of the diagonalization. 
 
 \begin{itemize}
 \item One can extract a sharp rate of convergence to a stationary  probability distribution, if there is one, see Diaconis et. al. \cite{dksc}. 
 \item The generator 
 $$G(z)=\left(\begin{array}{ccc}  &   &   \\  & G_{i,j}(z) &   \\  &   &  \end{array}\right)=F\left(\begin{array}{ccc}  &   &   \\  & -z \int_{-1}^1 {Q_i(\lambda) Q_j(\lambda) \over \lambda-z} d\psi(\lambda) &   \\  &   &  \end{array}\right)F^{T}$$
 \item One can use the Fokas, Its and Kitaev results, and benefit from the connection between orthogonal  polynomials and Riemann-Hilbert problems. 
 \item One can interpret random walks in random environment as a random spectral measure.
 \end{itemize}


\bibliographystyle{amsplain}

\end{document}